\documentclass[a4paper,11pt,twoside]{amsart}

\usepackage[latin1]{inputenc}
\usepackage{amsmath, amsfonts, amssymb,amsthm}
\usepackage[mathscr]{eucal} 
\usepackage[pdftex]{graphicx}

\usepackage[T1]{fontenc}
\usepackage[sc]{mathpazo}
\usepackage[all]{xy}

\usepackage{enumerate}

\newcommand{\N}{\mathbb{N}}
\newcommand{\Z}{\mathbb{Z}}
\newcommand{\R}{\mathbb{R}}
\newcommand{\C}{\mathbb{C}}
\newcommand{\Sp}{\mathbb{S}}
\newcommand{\h}{\mathbb{H}}

\newcommand{\Nil}{\mathrm{Nil}_3}

\newcommand{\pIm}{\mathrm{Im}}
\newcommand{\df}{\,\mathrm{d}}

\newcommand{\tq}{:\,}
\newcommand{\bnabla}{\nabla^b}

\newcommand{\prodesc}[2]{\left\langle #1, #2 \right \rangle}
\newcommand{\abs}[1]{\left\lvert #1 \right\rvert}

\newtheorem{theorem}{Theorem}
\newtheorem*{theorem*}{Theorem}
\newtheorem{proposition}{Proposition}
\newtheorem{corollary}{Corollary}
\newtheorem{lemma}{Lemma}

\theoremstyle{definition}
  \newtheorem{definition}{Definition}

\theoremstyle{remark}
  \newtheorem{remark}{Remark}

\numberwithin{equation}{section}

\title{Compact minimal surfaces in the Berger spheres}

\author{Francisco Torralbo}
\address{Departamento de Geometr\'{\i}a  y Topolog\'{\i}a \\
Universidad de Granada \\
18071 Granada, SPAIN} 
\email{ftorralbo@ugr.es}

\thanks{Research partially supported by a MCyT-Feder research project MTM2007-61775 and a Junta Andalucía  Grant P06-FQM-01642.}

\subjclass[2000]{Primary 53C42; Secondary 53C30}

\keywords{Surfaces, minimal, homogeneous 3-manifolds, Berger spheres}

\begin{document}
\maketitle

\begin{abstract}
We construct compact arbitrary Euler characteristic orientable and non-orientable minimal surfaces in the Berger spheres. Besides we show an interesting family of surfaces that are minimal in every Berger sphere, characterizing them by this property. Finally we construct, via the Daniel correspondence, new examples of constant mean curvature surfaces in $\Sp^2 \times \R$, $\h^2 \times \R$ and the Heisenberg group with many symmetries.
\end{abstract}

\section{Introduction}

The study of minimal surface in space forms is a classical topic in differential geometry. It is well know that there is no compact minimal surface in the Euclidean $3$-space and in the hyperbolic space. On the contrary, Lawson~\cite{Lawson70} showed that every compact surface but the projective plane can be minimally immersed in the $3$-sphere. His construction is based on the following fact: if a minimal surface contains a geodesic arc it can be continued as an analytic minimal surface by \emph{geodesic reflection} (see Definition~\ref{def:reflexion-geodesica}). Using this property he constructed compact minimal surfaces by reflecting the Plateau solution over a nice geodesic polygon. This method has been successfully applied to other spaces and other kinds of surface, like constant mean curvature ones (see, for example~\cite{KPS88, Karcher89, GB93}).

Besides the spaces form, the most regular $3$-manifolds are the homogeneous Riemannian ones. Between them, the compact type examples, which are expected to have compact minimal surfaces as in the space form case, are the Berger spheres. The author has been study the simplest constant mean curvature surfaces, in particular the minimal ones, in these family of homogeneous riemannian manifolds (cf.\ \cite{Tor10} and section~\ref{sec:new-examples}). Besides some studies about the stability of the constant mean curvature surfaces has been done (cf.\ \cite{TU10}).

In this paper we deal with the following problem: the construction of compact minimal surfaces in the Berger spheres. These $3$-manifolds are homogeneous Riemannian with isometry group of dimension $4$ which are, roughly speaking, $3$-spheres endowed with a family of deformation metrics of the round one (see section~\ref{sec:preliminaries} for details). The main trouble we are going to find is the following: \emph{the reflection over a geodesic of the Berger sphere could not be an ambient isometry, only a few of them are}.

First section is devoted to present in detail the family of Berger spheres and its geodesics and it shows what kind of geodesics satisfy that the reflection over them is an ambient isometry (cf.\ Lemma~\ref{lm:reflexion-isometria}). The second section characterizes a new family of minimal examples as the only ones that are minimal with respect to every Berger metric (see Proposition~\ref{prop:superficies-minimas-berger}). Section~\ref{sec:compact-orientable-examples} and \ref{sec:compact-non-orientable-examples} developed the main result (cf.\ Theorems~\ref{tm:compact-minimal-orientable}, \ref{tm:non-orientable} and Corollary~\ref{cor:compact-minimal-surfaces}) in the paper:
\begin{quote}
\itshape
Every compact surface but the projective plane can be minimally immersed in any Berger sphere. Furthermore, in the oriented case the immersion could be chosen without self intersection.
\end{quote}

Finally, in section~\ref{sec:sister-surfaces} we construct, as an application of the previous results and taking into account the Daniel correspondence~\cite{D}, complete constant mean curvature surfaces in the product spaces $\Sp^2 \times \R$ and $\h^2 \times \R$ (where $\Sp^2$ stands for the $2$-sphere and $\h^2$ for the hyperbolic plane) and in the Heisenberg group with many symmetries (see Proposition~\ref{prop:simetrias-correspondencia-Daniel}). In the product space case we are able to characterize which kind of symmetries these surfaces have (see Lemma~\ref{lm:simetrias-productos}).

This paper was prepared during the author's visit to IMPA (Brazil). The author wish to thank Harold Rosenberg for its valuable comments and suggestion and to Ildefonso Castro for the description of the minimal surfaces in the $3$-sphere invariant by a Killing field (see~\eqref{eq:integracion-ildefonso}).

\section{Preliminaries}\label{sec:preliminaries}
A Berger sphere is a usual $3$-sphere $\Sp^3 = \{(z, w) \in \C^2:\, \abs{z}^2 + \abs{w}^2 = 1\}$ endowed with the metric
\[
 g(X, Y) = \frac{4}{\kappa}\left[\prodesc{X}{Y} + \left(\frac{4\tau^2}{\kappa} - 1\right)\prodesc{X}{V}\prodesc{Y}{V} \right]
\]
where $\prodesc{\,}{\,}$ stands for the usual metric on the sphere, $V_{(z, w)} = J(z, w) = (iz, iw)$, for each $(z, w) \in \Sp^3$ and $\kappa$, $\tau$ are real numbers with $\kappa > 0$ and $\tau \neq 0$. For now on we will denote the Berger sphere $(\Sp^3,g)$ as $\Sp^3_b(\kappa, \tau)$. We note that if $\kappa = 4\tau^2$ then $\Sp^3_b(\kappa, \tau)$ is, up to homotheties, the round sphere. If $\kappa \neq 4\tau^2$ then the group of isometries of $\Sp^3_b(\kappa, \tau)$ is $\{A \in \mathrm{O}(4): AJ = \pm JA\}$ where $J$ is the matrix $J = \left(\begin{smallmatrix}J_1 & 0 \\ 0 & J_1\end{smallmatrix}\right)$ and  $J_1 = \left(\begin{smallmatrix}0 & -1 \\ 1 & 0\end{smallmatrix}\right)$. $\mathrm{Iso}(\Sp^3_b(\kappa, \tau))$ has two different connected components $\mathrm{U}_+(2)$ and $\mathrm{U}_-(2)$ given by
\[
\mathrm{U}_+(2) = \{A \in \mathrm{O}(4):\, AJ = JA\}, \quad
\mathrm{U}_-(2) = \{A \in \mathrm{O}(4):\, AJ = -JA\}
\]
and each of one is homeomorphic to the unitary group $\mathrm{U}(2)$. $\mathrm{U}_+(2)$, after identifying it with the unitary group, acts in the usual way over pair of complex numbers while $\mathrm{U}_-(2)$ acts by conjugation, that is, given $A = \left( \begin{smallmatrix} z_1 & z_2 \\ w_1 & w_2 \end{smallmatrix}\right) \in \mathrm{U}_-(2)$ then
\[
A(z, w) = (z_1 \bar{z} + z_2\bar{w}, w_1 \bar{z} + w_2 \bar{w})
\]

Taking into account the above relation between the standard metric on the sphere and the Berger metric, it is not difficult to get the next formula that links the Levi--Civita connection $\nabla$ of the round sphere to $\bnabla$, the one associated to the Berger metric:
\begin{equation}\label{eq:relacion-conexiones}
\bnabla_X Y = \nabla_X Y + \left(\frac{4\tau^2}{\kappa} - 1\right)[\prodesc{Y}{V}(JX)^\perp + \prodesc{X}{V}(JY)^\perp]
\end{equation}
where $J$ is the complex structure of $\C^2$, that is $J(z, w) = (iz, iw)$, and $()^\perp$ denotes the tangential component to the sphere.

The Hopf fibration $\Pi: \Sp^3_b(\kappa, \tau) \rightarrow \Sp^2(\kappa)$, where $\Sp^2(\kappa)$ stands for the $2$-sphere of radius $1/\sqrt{\kappa}$, 
\[
\Pi(z, w) = \frac{2}{\sqrt{\kappa}}\left(z\bar{w}, \frac{1}{2}(\abs{z}^2 - \abs{w}^2) \right),
\]  
is a Riemannian submersion whose fibers are geodesics. The vertical unit Killing field is given by $\xi = \frac{\kappa}{4\tau}V$. 

As the Berger spheres are homogeneous Riemannian manifolds to know all the geodesics it is enough to describe them at a fixed point, for example $(1,0) \in \Sp^3$. So, let $w = \left(i \frac{\kappa}{4\tau}\cos \theta, \frac{\sqrt{\kappa}}{2}\sin \theta e^{i\varphi} \right) \in \mathrm{T}_{(1,0)} \Sp^3_b$ an arbitrary unit vector, then the geodesic starting at $(1, 0)$ with speed $w$ is given by $\gamma_{\theta, \varphi}(s) = (z_1(s), z_2(s))$ where
\[
\begin{split}
z_1(s) &= \exp i\left(\frac{\kappa - 4\tau^2}{4\tau} \cos \theta\, s\right ) \left[\cos\left( \frac{\lambda \sqrt{\kappa}}{2} s\right) + i \frac{2\tau}{\sqrt{\kappa}} \frac{\cos \theta}{\lambda} \sin\left( \frac{\lambda \sqrt{\kappa}}{2} s\right)  \right] \\
z_2(s) &= \exp i\left( \varphi + \frac{\kappa - 4\tau^2}{4\tau}\cos \theta\, s\right) \frac{\sin \theta}{\lambda} \sin\left( \frac{\lambda \sqrt{\kappa}}{2} s\right)
\end{split}
\]
and $\lambda^2 = \sin^2 \theta + \frac{4\tau^2}{\kappa}\cos^2 \theta$ (cf.\ \cite{Rakotoniaina85}). We will say that a geodesic is \emph{horizontal} if its tangent vector is orthogonal to $V$ and \emph{vertical} if its tangent vector is co-linear with $V$. This terminology comes from the Hopf fibration, i.e., the vertical geodesics are the fibers of this Riemannian submersion. Notice that the horizontal and vertical geodesics are great circles. The ones starting at $(1, 0)$ are given by:
\begin{equation}\label{eq:parametrizacion-geodesicas}
\begin{split}
h_\varphi(s) &= \left(\cos \left(\frac{\sqrt{\kappa}}{2}s\right), \sin \left(\frac{\sqrt{\kappa}}{2}s\right) e^{i\varphi}\right), \qquad h'_\varphi(0) = \frac{\sqrt{\kappa}}{2}(0, e^{i\varphi}), \\
v(s) &= \left( e^{i\frac{\kappa}{4\tau}s}, 0\right), \qquad v'(0) = \frac{\kappa}{4\tau}(i, 0)
\end{split}
\end{equation}

\begin{remark}\label{rm:geodesics}
It is important to make the following remarks:
\begin{enumerate}
	\item Given two horizontal geodesics starting at $(1,0)$, $h_{\varphi_1}$ and $h_{\varphi_2}$, the angle between them is $\varphi_1 - \varphi_2$.
	\item The length of every vertical geodesic is $8\tau\pi/\kappa$, while the length of every horizontal geodesic is $4\pi/\sqrt{\kappa}$.
\end{enumerate}
\end{remark}

\begin{definition}\label{def:reflexion-geodesica}
A \emph{geodesic reflection} across a geodesic $\gamma$ of $\Sp^3_b(\kappa, \tau)$ is the map that sends a point $p$ to its ``opposite'' point on a geodesic through $p$ which meets $\gamma$ orthogonally. More precisely, if $\alpha$ is a geodesic that meet orthogonally $\gamma$ at $s = 0$ and $\alpha(s_0) = p$ then the geodesic reflection of $p$ with respect to $\gamma$ is the point $\alpha(-s_0)$.
\end{definition}

\begin{lemma}\label{lm:reflexion-isometria}
The geodesic reflection across an horizontal or a vertical geodesic is an isometry of $\Sp^3_b(\kappa, \tau)$.
\end{lemma}

\begin{proof}
As the Berger spheres are homogeneous Riemannian manifolds it is enough to check that the geodesic reflection across the vertical geodesic through $(1,0)$ and all the horizontal geodesic through $(1,0)$ are isometries.

In the first case it is easy to see that such transformation is given by $r_v(z, w) = (z,-w)$ which is clearly an isometry of the Berger sphere.

In the second case we fix first the horizontal geodesic given by $h_0$ (cf.\ \eqref{eq:parametrizacion-geodesicas}). In this special case it is not difficult to see that the geodesic reflection across $h_0$ is given by $R_0(z, w) = (\bar{z}, \bar{w})$ which is an isometry. Now, as the rotation of angle $\theta$ around the vertical geodesic $v$ through $(1,0)$ is $\rho(z, w) = (z, e^{i\theta}w)$ we can recover the geodesic reflection across any horizontal geodesic through $(1,0)$ by conjugation with this rotation. Therefore the geodesic reflection across the horizontal geodesic $h_\varphi$ (cf.\ \eqref{eq:parametrizacion-geodesicas}) is given by $R_\varphi(z, w) = (\bar{z}, e^{2i\varphi}\bar{w})$ which is an isometry. This finishes the proof.
\end{proof}

\section{New minimal surfaces in the Berger spheres}\label{sec:new-examples}

In \cite{Tor10} it is studied the simplest examples of constant mean curvature surfaces in the Berger sphere, indeed the \emph{rotationally} ones, that is, the constant mean curvature surfaces invariant by the $1$-parameter group of isometries given by $t\rightarrow \left(\begin{smallmatrix} 1 & 0 \\ 0 & e^{it} \end{smallmatrix}\right)$. Among them it was found the following minimal ones:
\begin{itemize}
	\item The equator $\{(z, w) \in \Sp^3_b(\kappa, \tau)\tq \pIm(z) = 0\}$, which is the only minimal sphere, up to ambient isometries, in $\Sp^3_b(\kappa, \tau)$,
	\item the Clifford torus $\{(z, w) \in \Sp^3_b(\kappa, \tau):\abs{z}^2 = \abs{w}^2 = \frac{1}{2}\}$,
	\item a $1$-parameter family of examples, call \emph{unduloids} because they looks like euclidean unduloids. Some of them produces tori and, even more, for some Berger spheres with small bundle curvature $\tau$ (with respect to a fixed $\kappa$) there exist embedded ones (cf.\ \cite[Remark 3.2]{Tor10}).
\end{itemize}

The first two examples are in fact minimal surfaces in every Berger sphere, that is, they are minimal surfaces in $\Sp^3_b(\kappa, \tau)$ for every $\kappa > 0$ and $\tau \neq 0$. So a natural question arises: \emph{Does exist more surfaces that are simultaneously minimal  in all the Berger spheres?} This section is devoted to answer positively this question providing a new $1$-parameter family of minimal surfaces in the Berger spheres (see Proposition~\ref{prop:superficies-minimas-berger}).

Firstly we are going to related the mean curvature of an arbitrary surface in a Berger sphere with its mean curvature as a surface in the usual $3$-sphere.

\begin{lemma}\label{lm:relacion-curvaturas-medias}
Let $\Phi: \Sigma \rightarrow \Sp^3$ be an immersion of a surface $\Sigma$. Then the mean curvatures $H$ and $H^b$ of $\Phi$ in the round sphere $(\Sp^3,\prodesc{\,}{\,})$ and the Berger sphere $\Sp^3_b(\kappa, \tau)$ respectively are related by
  \[
    H^b = \frac{\sqrt{\kappa}}{2} \frac{1}{\sqrt{1 - \left(1 - \frac{\kappa}{4\tau^2} \right)\nu^2}} \left[H + \frac{\left(1 - \frac{\kappa}{4\tau^2}\right)\prodesc{\nabla \nu}{V}}{2\left[1 - \left(1 - \frac{\kappa}{4\tau^2} \right)\nu^2 \right]} \right]
  \]
  where $\nu = \prodesc{N}{V}$, $N$ is a unit normal vector field of $\Sigma$ in $(\Sp^3, \prodesc{\,}{\,})$ and $V_{(z,w)} = (iz, iw)$.
\end{lemma}

\begin{proof}
First if $N$ and $N^b$ are two unit normal vector field of $\Phi$ in $\Sp^3$ and $\Sp^3_b(\kappa, \tau)$ respectively then they are related by the equality:
  \[
    N^b = \frac{\sqrt{\kappa}}{2\sqrt{1 - \left(1 - \frac{\kappa}{4\tau^2}\right) \nu^2}} \bigl[N - \nu\left(1 - \frac{\kappa}{4\tau^2} \right)V \bigr],
  \]
where $\nu = \prodesc{N}{V}$.

Now taking into account~\eqref{eq:relacion-conexiones} and the above relation between the normal vector fields we get
  \[
    \begin{split}
    \frac{\kappa}{4f} g(\bnabla_{e_i} N^b, e_j) = &\left(\frac{4\tau^2}{\kappa} - 1 \right)\prodesc{e_i}{V}\prodesc{e_j}{JN} +\\ 
    &-\lambda_i\left[\delta_{ij} + \left(\frac{4\tau^2}{\kappa} - 1 \right) \prodesc{e_i}{V}\prodesc{e_j}{V} \right] + \\
    & -\left(\frac{4\tau^2}{\kappa} - 1 \right)\prodesc{e_j}{V}e_i(\nu),
    \end{split}
  \]
  where $\{e_1, e_2\}$ is an orthonormal reference with respect to $\prodesc{\,}{\,}$ that diagonalize the shape operator of $\Phi:\Sigma \rightarrow (\Sp^3, \prodesc{\,}{\,})$, $\lambda_i$ are the principal curvatures and $f:\Sigma \rightarrow \R$ is the function given by $f = \sqrt{\kappa}\left/2\sqrt{1-\left(1 - \frac{\kappa}{4\tau^2} \right)\nu^2}\right.$.
  Finally, as $-2H^b = \sum_{i,j = 1}^2 g^{ij} g(\bnabla_{e_i} N^b, e_j)$ a straightforward computation yields the result.
\end{proof}

Thanks to the above result if $\Sigma$ is a minimal surface in the round sphere then it is minimal in some Berger sphere $\Sp^3_b(\kappa, \tau)$ (in fact, in every Berger sphere) if and only if $\prodesc{\nabla \nu}{V} = 0$. The equator $\{(z, w) \in \Sp^3:\, \pIm(z) = 0\}$ and the Clifford torus $\{(z, w) \in \Sp^3:\, \abs{z}^2  =\abs{w}^2 = \frac{1}{2}\}$ satisfy this condition so they are, as we already knew, minimal surfaces in every Berger sphere. Next proposition classifies all the surfaces that are simultaneously minimal in every Berger sphere.

\begin{proposition}\label{prop:superficies-minimas-berger}
Let $\Phi: \Sigma \rightarrow \Sp^3$ be an immersion of a surface $\Sigma$. Suppose that $\Phi$ is a minimal immersion with respect to any Berger metric. Then $\Phi$ is congruent to a piece of one of the immersion $\Phi_{c}: \R^2 \rightarrow \Sp^3$ given by:
	\[
	\Phi_{c}(s, t) = \Bigl(\cos(s) e^{ict},\, \sin(s) e^{it} \Bigr),
	\]
	where $c \in \R$.
\end{proposition}

\begin{remark}~\label{rm:superficies-minimas-berger}
\begin{enumerate}
	\item If $c = 0$ we get the equator $\{(z, w) \in \C^2:\pIm(z) = 0\}$.

	\item The immersions $\Phi_{c}$ are invariant by the $1$-parameter group of isometries of $U(2)$ given by:
	\[
	\left\{
		\begin{pmatrix}
			e^{ict}	&	0	\\
			0		&	e^{it}
		\end{pmatrix}
		\tq t \in \R
	\right\}.
	\]
	The immersions $\Phi_{c}$ do not coincide with the family of minimal unduloids described in \cite[Theorem 1.(iii)]{Tor10} because they are not invariant by the group of isometries considered there, except for the sphere where $c = 0$ (see previous item).
	
	\item If $c = m/n$ with $m,n \in \Z$ we get the compact examples $\tau_{m,n}$ constructed by Lawson~\cite[Section 7, p.\ 350]{Lawson70}. It is interested to remark that (cf.~\cite[Theorem 3]{Lawson70} for more details):
	\begin{itemize}
		\item To each pair of coprime natural numbers $m, n$,  $\tau_{m,n}$ is a compact surface with zero Euler characteristic. 
		\item $\tau_{m, n}$ is a Klein bottle if and only if $m\cdot n$ is an even number.
		\item $\tau_{1,1}$ is the Clifford torus and it is the only surface of this type that is embedded. More precisely, up to an isometry of $U(2)$, we get
	\[
	\frac{1}{\sqrt{2}} 
	\begin{pmatrix}
		1	&	i 	\\
		1	&	-i
	\end{pmatrix}
	\Phi_{1}(s, t) = \frac{1}{\sqrt{2}}\Bigl(e^{i(t+s)}, e^{i(t-s)} \Bigr)
	\]
	\end{itemize}
\end{enumerate}
\end{remark}

\begin{proof}
Let $\Sigma$ be a minimal surface in every Berger sphere. Then, thanks to Lemma~\ref{lm:relacion-curvaturas-medias}, it must be $\prodesc{\nabla \nu}{V} = 0$. Let $z = x + iy$ a conformal parameter for $\Phi: \Sigma \rightarrow (\Sp^3, \prodesc{\,}{\,})$ with conformal factor $e^{2u}$. Then, the compatibility equations (cf.\ \cite[Theorem 2.3, p. 283]{FM}) of the immersion are given, for $\kappa = 4$, $\tau = 1$ (which is the case of the round sphere) and $H = 0$, by:
\begin{equation}\label{eq:compatibilidad}
\begin{aligned}
  p_{\bar{z}} &= 0, &  \nu_z &= iA - 2e^{-2u} \bar{A}p \\
  A_{\bar{z}} &= i\frac{e^{2u}}{2}\nu,\quad  &
  \abs{A}^2 &= \frac{e^{2u}}{4}(1 - \nu^2)
\end{aligned}
\end{equation}
where $A = \prodesc{\Phi_z}{V}$ and $p = \prodesc{\Phi_{zz}}{N}$. As the Hopf differential $\Theta = p(\mathrm{d}z)^2$ is holomorphic we can normalized it by $\lambda i$, with $\lambda \in \R$. Now, the condition $\prodesc{\nabla \nu}{V} = 0$ becomes in
\[
0  = \nu_z \bar{A} + \nu_{\bar{z}}A = -4\lambda e^{-2u}\pIm(A^2)
\]
If $\lambda = 0$ then $\Sigma$ is a totally umbilical minimal surface so it must be a totally geodesic sphere which correspond to the case (1) in Remark~\ref{rm:superficies-minimas-berger}.

We suppose now that $\lambda \neq 0$. The last equation allows us to suppose that the real part of $A$ is zero. So, using~\eqref{eq:compatibilidad} it must be
\[
A = i\frac{e^{u}}{2}\sqrt{1-\nu^2}
\]
In that case from the equation for $\nu_z$ in~\eqref{eq:compatibilidad} we get
\[
\nu_z = -\frac{\sqrt{1-\nu^2}}{2}(e^u + 2\lambda e^{-u})
\]
From the above equation we deduce that $\nu_y = 0$, i.e., $\nu = \nu(x)$. Using again the last equation we get that $u = u(x)$ too.

It is well known that the minimal surfaces of the round sphere are controlled by the sinh-Gordon equation. In this case, as the conformal parameter depends only on one variable they must satisfy the equation
\begin{equation}\label{eq:sinh-Gordon}
u'' + (e^{2u} - 4\lambda^2e^{-2u}) = 0, 
\end{equation}
whose first integral is given by
\begin{equation}\label{eq:sinh-Gordon-energia}
(u')^2 + (e^{2u} - 4\lambda^2e^{-2u}) = E
\end{equation}

The above equation has been intensively studied in the literature. Its solution can be described in terms of Jacobi elliptic functions. In fact, it is possible to integrate explicitly the minimal immersions that they produce as follows:
\begin{equation}\label{eq:integracion-ildefonso}
\Psi_{a, b} = \frac{1}{\sqrt{a^2 - b^2}}\left(\sqrt{e^{2u(x)}-b^2} e^{i(bF(x) + ay)},\sqrt{a^2 - e^{2u(x)}}e^{i(aG(x) +by)}\right)
\end{equation}
where
\[
\begin{split}
F(x) &= \int_0^x \frac{\sqrt{P(e^{2u(t)})}}{e^{2u(t)} - b^2} \df t, \quad 
G(x) = \int_0^x \frac{\sqrt{P(e^{2u(t)})}}{a^2 - e^{2u(t)}} \df t, \\
P(e^{2u(t)}) &= (a^2 - e^{2u(t)})(e^{2u(t)} - b^2) - e^{2u(t)}u'(t)^2 
\end{split}
\]
and $u(x)$ is a solution of~\eqref{eq:sinh-Gordon} for $\lambda^2 = a^2b^2/4$ and $E = a^2 + b^2$.

Now, among all these solutions the only ones that satisfy the extra condition that $\prodesc{\nabla \nu}{V} = 0$ are precisely, up to reparametrization, the family $\Phi_{c}$ and the proof finishes.
\end{proof}

\section{Compact minimal surfaces with arbitrary topology}\label{sec:compact-orientable-examples}
In this section we are going to construct new examples of compact minimal surfaces in the Berger sphere following a technique developed by Lawson in~\cite{Lawson70}. The main idea is to solve the Plateau problem over a nice geodesic polygon and to successively reflect the obtained minimal surface over the border geodesics to generated a compact minimal surface.

The key result that we are going to apply is the following \emph{reflection principle}:
\begin{quote}
\itshape
Given a minimal surface $\Sigma$ on a Berger sphere $\Sp^3_b(\kappa, \tau)$, if $\partial\Sigma$ contains a horizontal or a vertical geodesic arc $\gamma$ then $\Sigma$ can be continued as an analytic minimal surface across each non-trivial component of $\partial\Sigma \cap \gamma$ by geodesic reflection across $\gamma$. 
\end{quote}

In deep, in the above situation if we denote by $r$ the geodesic reflection across the border geodesic $\gamma$ then both $\Sigma$ and $r(\Sigma)$ are minimal surfaces because $r$ is an isometry (cf.\ Lemma~\ref{lm:reflexion-isometria}). Even more, their tangent planes coincide along $\gamma$ because the map $r$ change the co-normal vector field $\eta$ along $\gamma$ to $-\eta$. Therefore, as they are minimal surfaces, $\Sigma$ can be continue as an analytic minimal surface by $r(\Sigma)$.

\subsection{Choosing the polygon}
\begin{figure}[htbp]
\centering
\includegraphics{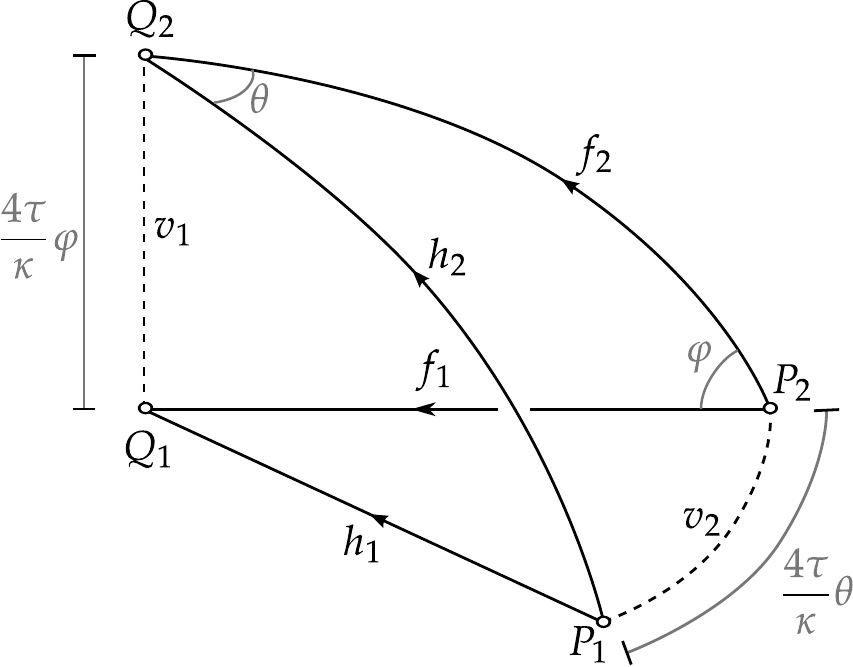}
\caption{Polygon $\Gamma$ by stereographic projection from the point $(0, -1)$\label{fig:polygon}}
\end{figure}

To help us to visualize the polygon we consider the stereographic projection of the sphere from the point $(0,-1)$. The vertical geodesic $v_1$ from $(0, 1)$ becomes the $z$ axis and the one from $(1,0)$, $v_2$, becomes the unit circle centered at the origin in the $xy$ plane. Choose $P_1, P_2$ two points in $v_2$ such that $\mathrm{distance}(P_1,P_2) = (4\tau/\kappa) \theta$ and $Q_1$, $Q_2 \in v_1$ such that $\mathrm{distance}(Q_1,Q_2) = (4\tau/\kappa) \varphi$ (see figure~\ref{fig:polygon}).

We are going to consider the polygon $\Gamma = Q_1P_1Q_2P_2$ where its edges are the horizontal geodesics $f_1$, $f_2$, $h_1$ and $h_2$ given by:
\begin{align*}
f_1(t) &= (\cos(t) e^{i\theta}, \sin(t)), & f_2(t) &= (\cos(t)e^{i\theta}, \sin(t) e^{i\varphi})	& t \in [0,\pi/2]\\
h_1(t) &= (\cos(t), \sin(t)), & h_2(t) &= (\cos(t), \sin(t) e^{i\varphi})	& t \in [0,\pi/2]
\end{align*}

\begin{remark}\label{rmk:rotation-double-angle}
It is not difficult to check that the group $G_\Gamma$ generated by geodesic reflections around the edges of $\Gamma$ is generated by the maps $H_1(z, w) = (\bar{z}, \bar{w})$, $H_2(z, w) = (\bar{z}, e^{2i\varphi}\bar{w})$, $F_1(z, w) \mapsto (e^{2i\theta}\bar{z}, \bar{w})$ and $F_2(z, w) = (e^{2i\theta}\bar{z}, e^{2i\varphi}\bar{w})$.

Moreover, two successively reflexions around the vertices $P_1$ or $P_2$ (resp.  $Q_2$ or $Q_2$) produce a rotation of angle $2\varphi$ (resp. $2\theta$) around the vertical geodesic $v_2$ (resp. $v_1$), i.e., the  maps $(z, w) \mapsto (z, e^{2i\varphi}w)$ and $(z,w) \mapsto (e^{2i\theta}z, w)$ are in the group $G_\Gamma$.
\end{remark}

\subsection{Solving the Plateau problem}
In order to ensure the existence of a minimal surface with border $\Gamma$ we are going to use the following general result due to Meeks and Yau~\cite{MY80}.

\begin{theorem*}
Let $W$ be a compact Riemannian $3$-manifold with piecewise smooth mean convex boundary\footnote{Suppose that W is a compact Riemannian $3$-manifold that embeds in
the interior of another Riemannian $3$-manifold. W is said to have \emph{piecewise smooth mean convex boundary} if $\partial W$ is a two-dimensional complex consisting of a finite number of smooth two-dimensional compact simplices with interior angles less than or equal to $\pi$, each one with non negative mean curvature with respect to the inward pointing normal}. Let $\Gamma$ be a smooth collection of pairwise disjoint closed curves in $\partial W$, which bounds a compact orientable surface in $W$. Then, there exists an embedded orientable surface $\Sigma \subseteq W$ with $\partial \Sigma = \Gamma$ that minimizes area among all orientable surfaces with the same boundary.
\end{theorem*}

Consider $W$ the solid \emph{tetrahedron} determined by $\Gamma$, $v_1$ and $v_2$ (see figure~\ref{fig:polygon}). Clearly $W$ is a compact Riemannian $3$-manifold whose border is a two-dimensional complex consisting of four two-dimensional compact simplices, which are going to call \emph{faces}, with interior angles less than or equal to $\pi$ if and only if $\varphi, \theta \leq \pi$. If we denote every face of $\partial W$ by $\blacktriangle ABC$ where $A$, $B$ and $C$ are the vertex, then we can parametrized them by the following:
\begin{align*}
	\Phi_1: [0, \frac{\pi}{2}]\times[0, \varphi] &\rightarrow \blacktriangle Q_1Q_2P_1 \subseteq \Sp^3_b(\kappa, \tau), & \Phi_1(t, s) &= \bigl( \cos(t), \sin(t)e^{is}\bigr) \\
	\Phi_2: [0, \frac{\pi}{2}]\times[0, \varphi] &\rightarrow \blacktriangle Q_1Q_2P_2 \subseteq \Sp^3_b(\kappa, \tau), & \Phi_2(t, s) &= \bigl( \cos(t)e^{i\theta}, \sin(t)e^{is}\bigr) \\
	\Phi_3: [0, \frac{\pi}{2}]\times[0, \theta] &\rightarrow \blacktriangle Q_1P_1P_2 \subseteq \Sp^3_b(\kappa, \tau), & \Phi_3(t, s) &= \bigl( \cos(t)e^{is}, \sin(t)\bigr) \\
	\Phi_4: [0, \frac{\pi}{2}]\times[0, \theta] &\rightarrow \blacktriangle Q_2P_1P_2 \subseteq \Sp^3_b(\kappa, \tau), & \Phi_4(t, s) &= \bigl( \cos(t)e^{is}, \sin(t)e^{i\varphi}\bigr)
\end{align*}
Then the only condition that we have to check in order to $W$ has a \emph{piecewise smooth mean convex boundary} is that every face has non-negative mean curvature with respect to the inward pointing normal.

First of all $\Phi_1$ is a minimal immersion because is a small piece of the minimal sphere $\{(z, w) \in \Sp^3_b(\kappa, \tau):\, \pIm(z) = 0\}$. On the other hand $\Phi_2 = \rho_\theta \circ \Phi_1$, where $\rho_\theta(z, w) = (e^{i\theta}z, w)$ is the rotation around the vertical geodesic $v_1$ of angle $\theta$ so $\Phi_2$ is minimal too. Moreover $\Phi_3 = s\circ \Phi_1$, where $s(z, w) = (w, z)$ is an isometry of the Berger sphere and $\Phi_4 = \rho_\varphi \circ \Phi_3$, where $\rho_\varphi(z, w) = (z, e^{i\varphi}w)$ is the rotation of angle $\varphi$ around the  vertical geodesic $v_2$. Then both $\Phi_3$ and $\Phi_4$ are minimal immersions too. To sum up $W$ has piecewise smooth mean convex boundary so the above theorem ensures us that there exist a minimal surface $\Sigma$ with border $\Gamma$ if and only if $\varphi, \theta \leq \pi$.

\subsection{Generating the surface}
In the above situation take $n, m \in \N$ and choose $\varphi = \pi/(n+1)$ and $\theta = \pi/(m+1)$. We denote by $\Sigma_\Gamma^{(m,n)}$ a minimal surface in $\Sp^3_b(\kappa, \tau)$ with border $\Gamma$. Let $G_\Gamma$ the group of isometries generated by geodesic reflection around the edges of $\Gamma$. Then, thanks to the choice of $\varphi$ and $\theta$, $G_\Gamma$ is a discrete group so $\Sigma^{(m,n)} = \bigcup_{g \in G_\Gamma} g(\Sigma_\Gamma^{(m,n)})$ is a compact minimal surface in $\Sp^3_b(\kappa, \tau)$. Even more, if $H = \{g \in G_\Gamma:\, g(\Sigma_\Gamma^{(m,n)}) = \Sigma_\Gamma^{(m,n)}\}$, then $\mathrm{ord}(G)/\mathrm{ord}(H) = 2(m+1)(n+1)$ because $2(n+1)$ reflections at $v_2$ produce the identity and $(m+1)$ reflections at $v_1$ of the resulting surface produce the whole surface $\Sigma^{(m,n)}$.

Moreover if we denote by $K$ the Gauss curvature of $\Sigma_\Gamma^{(m,n)}$ then it is easy to check that, by the Gauss-Bonnet Theorem,
\[
\int_{\Sigma_\Gamma^{(m,n)}} K = \frac{2\pi}{(m+1)(n+1)}(1-mn)
\]
and so
\[
2\pi \chi(\Sigma^{(m,n)}) = \int_{\Sigma^{(m,n)}} K = \frac{\mathrm{ord}(G)}{\mathrm{ord}(H)}\int_{\Sigma_\Gamma^{(m,n)}}K = 4\pi(1-mn)
\]
where $\chi(\Sigma^{(m,n)})$ denotes the Euler characteristic. Finally we get that the surface $\Sigma^{(m,n)}$ is oriented and its genus is $mn$. To sum up we state the following result.

\begin{theorem}\label{tm:compact-minimal-orientable}
For every $g \geq 0$ there exist a compact embedded minimal surface of genus $g$ in $\Sp^3_b(\kappa, \tau)$. 
\end{theorem}

\begin{proof}
By the construction method we already know that there exists a compact minimal surface $\Sigma^{(m,n)}$ (possibly with singularities) in $\Sp^3_b(\kappa, \tau)$ of genus $g = mn$, $m,n \in \N$. 

The embeddedness follows from the fact that the fundamental piece $\Sigma^{(m,n)}_\Gamma$ fits in the \emph{solid tetrahedron} determined by the border $\Gamma$ and this tetrahedron does not never intersect itself as we go on reflecting it by the isometries of the group $G_\Gamma$.

Finally we prove that no singularity appears. We know that $\Sigma^{(m,n)}$ is regular up to the border $\Gamma$ and at the interior of the edges of $\Gamma$ the surface is smooth too so we only have to check the regularity at the vertices $P_1$, $P_2$, $Q_1$, $Q_2$ and their reflections. But this is a consequence of a more general \emph{removable singularity result} (cf.\ \cite[Proposition 1]{CS85}). In order to apply that result we only notice that the surface $\Sigma^{(m,n)}$ is locally embedded around the vertex (it is, in fact, embedded).
\end{proof}

\begin{remark}\label{rm:ejemplos-lawson-genero-1}~ 
\begin{enumerate}
	\item We can produce other type of compact minimal surface by just changing the geodesic polygon over the tetrahedron $\partial W$. In deep, take the geodesic polygon $\tilde{\Gamma}$ which join the points $P_1Q_1Q_2P_2$ in figure~\ref{fig:polygon} with angles $\varphi = \pi/(n+1)$ and $ \theta = \pi/(m+1)$. Then, by the previous process, we can produce a compact minimal surface that, in this case, always has Euler characteristic zero. In fact the resulting surface can be explicitly described by the immersion $\Phi_{m/n}$ given in Proposition~\ref{prop:superficies-minimas-berger}.

	\item Let $r \in \N$ and $T_r$ the group generated by the map 
\[
(z, w) \rightarrow (e^{\frac{2\pi i}{r}}z, e^{\frac{2\pi i}{r}} w).
\]
Notice that $T \subset \mathrm{Iso}(\Sp^3_b(\kappa, \tau))$ for every $\kappa$ and $\tau$. We will denote  the quotient $\Sp^3_b(\kappa, \tau)/T_r$ with the induced metric by $L^r_b(\kappa, \tau)$. Then $L^r_b(\kappa, \tau)$ is a homogeneous riemannian\footnote{The isometry $\left(\begin{smallmatrix}a & -\bar{b} \\ b & \bar{a}\end{smallmatrix}\right) \left( \begin{smallmatrix} \bar{\alpha} & \bar{\beta} \\ - \beta & \alpha \end{smallmatrix}\right)$ sends the point $(\alpha, \beta)$ to $(a, b)$ and it is induced to the quotient $L^r_b(\kappa, \tau)$ because it commutes with $T_r$} $3$-manifold with isometry group of dimension $4$.

For each $r\in \N$ such that $r$ divides $n+1$ and $m+1$ the surface $\Sigma^{(m,n)}$ is invariant by the group $T_r$ an so it is induced to the quotient $L_b^r(\kappa, \tau)$ as a compact minimal surface with genus $1 - (1-mn)/r$.
\end{enumerate}
\end{remark}

\section{Non-orientable compact minimal surfaces}\label{sec:compact-non-orientable-examples}
\begin{figure}[htbp]
\centering
\includegraphics[width=0.9\textwidth]{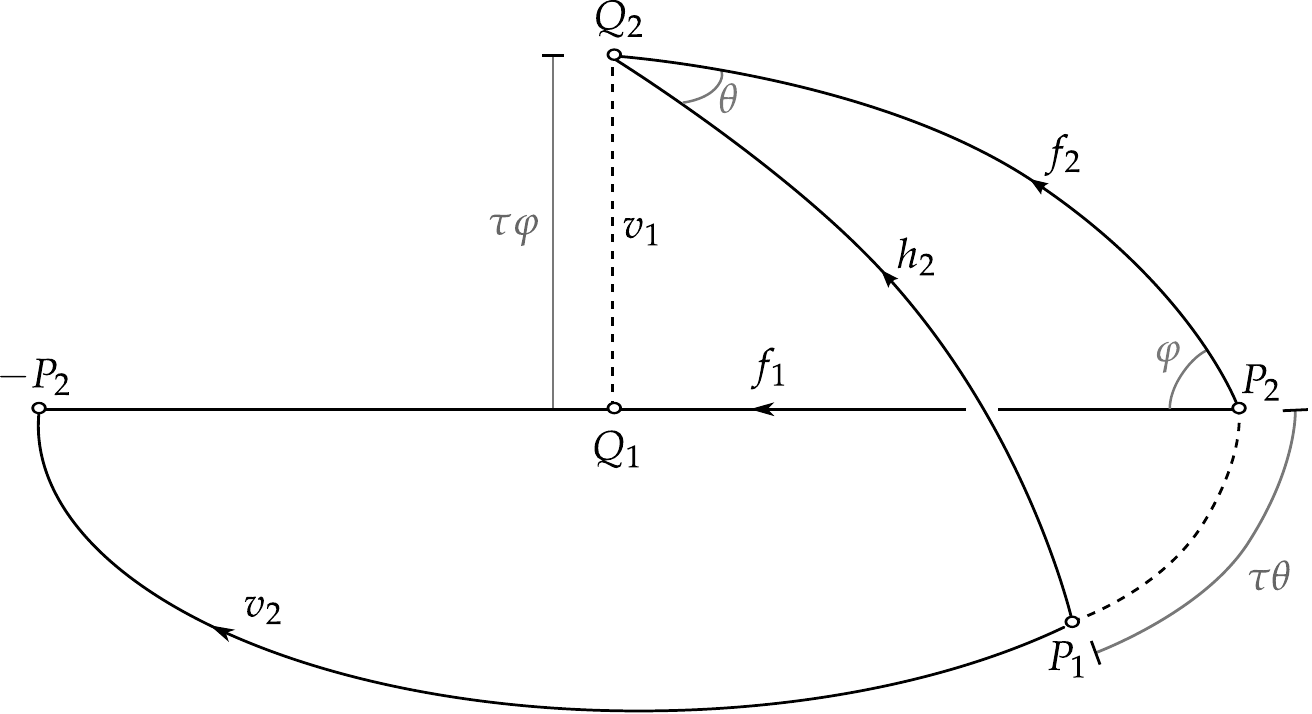}
\caption{Polygon $\Gamma'$ by stereographic projection from the point $(0, -1)$\label{fig:polygon2}}
\end{figure}

Now we are going to construct, following the same procedure as in the previous section, compact non-orientable minimal surfaces. We start with the points $P_1$, $P_1$, $Q_1$ and $Q_2$ as in section 4 and we define a polygon $\Gamma' = P_1Q_2P_2(-P_2)$ where its edges are the horizontal geodesics $f_1$, $f_2$, $h_2$ and the vertical one $v_2$ parametrized by (see figure~\ref{fig:polygon2}):
\begin{align*}
f_1(t) &= (\cos(2t) e^{i\theta}, \sin(2t)), & f_2(t) &= (\cos(t)e^{i\theta}, \sin(t) e^{i\varphi})	& t \in [0,\pi/2]\\
h_2(t) &= (\cos(t), \sin(t) e^{i\varphi})	& & & t \in [0,\pi/2] \\
v_2(t) &= (e^{it}, 0) & & & t \in [-\theta, 0]
\end{align*}

Consider $W$ the region delimited by the minimal spheres $\{(z, w) \in \Sp^3:\, \pIm(w) = 0\}$, $\{(z, w) \in \Sp^3:\, \pIm(e^{-i\theta}z) = 0\}$ and $\{(z, w) \in \Sp^3:\, \pIm(e^{-i\varphi}w) = 0\}$. Then the polygon $\Gamma'$ is in the border of $W$ and $W$ has piecewise smooth mean convex boundary if and only if the angles $\varphi, \theta \leq \pi$. 

Take $n,m \in \N$ with $n$ \textbf{odd} and choose $\varphi = \pi/(n+1)$ and $\theta = \pi/(m+1)$. Therefore there exist an embedded minimal surface $\Lambda^{(m,n)}_{\Gamma'}$ with border $\Gamma'$. Let $G_{\Gamma'}$ the group of isometries generated by geodesic reflections around the edges of $\Gamma'$. Then, thanks to the choice of $\varphi$ and $\theta$, $G_{\Gamma'}$ is a discrete group so $\Lambda^{(m,n)} = \bigcup_{g \in G_{\Gamma'}} g(\Lambda^{(m,n)}_{\Gamma'})$ is a compact minimal surface (possibly with singularities) in $\Sp^3_b(\kappa, \tau)$. 

Let $\Lambda$ by the set constructed by reflecting $\Lambda^{(m,n)}_{\Gamma'}$ $2(m+1)$ times at $Q_2$ and then reflecting the resulting configuration $(n+1)$ times at $P_2$. We claim that $\Lambda = \Lambda^{(m,n)}$. It is sufficient to prove that $\Lambda$ is invariant under the reflection around the edges of $\Gamma'$. Firstly it is clear by the construction of $\Lambda$ that it is invariant under the reflection under the edges $Q_2P_2$, $Q_2P_1$ and $Q_1P_2$. Furthermore, as $n$ is odd, the reflection around $P_1(-P_2)$ coincide with the $(n+1)/2$ times reflection around $P_2$ so the claim is proved.

By the previous reasoning we also know that if $H = \{g \in G_{\Gamma'}:\, g(\Lambda^{(m,n)}_{\Gamma'}) = \Lambda^{(m,n)}_{\Gamma'}\}$, then $\mathrm{ord}(G)/\mathrm{ord}(H) = 2(m+1)(n+1)$.

Finally, if $K$ denote the Gaussian curvature of $\Lambda^{(m,n)}_{\Gamma'}$ then, by the Gauss-Bonnet theorem
\[
\int_{\Lambda^{(m,n)}_{\Gamma'}} K = \frac{\pi}{(n+1)(m+1)}(1-mn)
\]
and so
\[
2\pi \chi(\Lambda^{(m,n)}) = \int_{\Lambda^{(m,n)}} K = \frac{\mathrm{ord}(G)}{\mathrm{ord}(H)}\int_{\Lambda^{(m,n)}_{\Gamma'}} K = 2\pi(1 - mn)
\]
That is, $\chi(\Lambda^{(m,n)}) = 1 - mn$. It is clear that if $1- mn$ is odd the resulting surface is non-orientable. Moreover, if $1 -mn$ is even the surface is also non-orientable. To prove this we use the same argument as Lawson: we consider a non-contractible, closed curve on the piece of the surface $\Lambda^{(m,n)}_{\Gamma'} \cup R(\Lambda^{(m,n)}_{\Gamma'}) \cup U$ where $R$ is the reflection around $P_1(-P_2)$ and $U$ is a small neighborhood of $P_2$ on the surface $\Lambda^{(m,n)}$. Such a curve is orientation reversing.

\begin{theorem}\label{tm:non-orientable}
For every pair of natural numbers $n$ and $m$, with $n$ odd, there exists a non-orientable, compact minimal surfaces $\Lambda^{(m,n)}$ with Euler characteristic $1-mn$
\end{theorem}

\begin{proof}
The previous construction ensures that there exist a compact (possibly with singularities) non-orientable minimal surface $\Lambda^{(m,n)}$ with Euler characteristic $1-mn$. Now, if we look locally at a vertex, the surface only go once around it as we go reflecting the small piece $\Lambda^{(m,n)}_{\Gamma'}$ by the edges of $\Gamma'$ so $\Lambda^{(m,n)}$ is locally embedded around the vertex. Hence, by the same argument as in the proof of Theorem~\ref{tm:compact-minimal-orientable} we know that the surface $\Lambda^{(m,n)}$ is smooth. Finally no branching points appear at the vertices because, as we notice before, the surface is locally embedded around them.
\end{proof}

As a corollary from Theorem~\ref{tm:compact-minimal-orientable} and Theorem~\ref{tm:non-orientable} we get the following:

\begin{corollary}\label{cor:compact-minimal-surfaces}
Every compact surface but the projective plane can be minimally immersed into any Berger sphere $\Sp^3_b(\kappa, \tau)$. Even more, for orientable surfaces the immersions can be chosen without self-intersections.
\end{corollary}

\begin{remark}
The projective plane is forbidden because the only minimal sphere, up to congruences, is the equator $\{(z, w):\, \pIm(z) = 0\}$ which is embedded. 
\end{remark}

\section{Constant mean curvature surfaces in $\Sp^2 \times \R$, $\h^2 \times \R$ and the Heisenberg group}\label{sec:sister-surfaces}

An important tool in the description of constant mean curvature surfaces in space forms is the classical Lawson correspondence. It establishes an isometric one-to-one local correspondence between constant mean curvature surfaces in different space forms that allows to pass, for instance, from
minimal surfaces in $\Sp^3$ to constant mean curvature $1$ surfaces in $\R^3$. Lawson used it to construct double periodic constant mean curvature surfaces in $\R^3$ from the compact minimal ones that he constructed in $\Sp^3$ (cf.\ \cite[Section 14]{Lawson70}).

In 2007 Daniel~\cite[Theorem 5.2]{D} generalized the Lawson correspondence to the context of homogeneous Riemannian spaces. The author related constant mean curvature immersions between the simply connected homogeneous Riemannian spaces with isometry group of dimension $4$, and he called these immersions \emph{sister immersions}.

Although the Daniel correspondence is more general we are going to consider only the isometric correspondence between minimal surfaces in the Berger spheres and constant mean curvature surfaces in $\Sp^2 \times \R$, $\h^2 \times \R$, where $\Sp^2$ and $\h^2$ stand for the $2$-sphere and the hyperbolic plane of constant curvature $1$ and $-1$ respectively, and the Heisenberg group $\Nil$. We can state it as follows: 

\begin{quote}
\itshape
Given a simply connected surface $\Sigma$ then
\begin{itemize}
	\item To each minimal immersion of $\Sigma$ in $\Sp^3_b(4\tau^2 + 1, \tau)$ it corresponds a constant mean curvature $H = \tau$ isometric immersion of $\Sigma$ in $\Sp^2 \times \R$,
	\item to each minimal immersion of $\Sigma$ in $\Sp^3_b(4\tau^2 - 1, \tau)$ ($\tau^2$ must be greater than $1/4$) it corresponds a constant mean curvature $H = \tau$ isometric immersion of $\Sigma$ in $\h^2 \times \R$.

	\item To each minimal immersion $\Sigma$ in $\Sp^3_b(\kappa, \tau)$ with $\kappa - 4\tau^2 < 0$ it corresponds a constant mean curvature $H = \sqrt{\kappa}/2\sqrt{4\tau^2 - \kappa}$ isometric immersion of $\Sigma$ in $\Nil$.
\end{itemize}
\end{quote}

Hence, the minimal immersions that we constructed in Theorem~\ref{tm:compact-minimal-orientable} in the family of Berger spheres produce, passing through the universal cover, complete constant mean curvature immersions in $\Sp^2 \times \R$, $\h^2 \times \R$ and $\Nil$ with arbitrary mean curvature. As the initial surface has many symmetries the corresponding sister one should have many too. The following proposition ensures it.

\begin{proposition}\label{prop:simetrias-correspondencia-Daniel}
Let $\Sigma^{(m,n)}$ a compact minimal surface of $\Sp^3_b(\kappa, \tau)$ given in Theorem~\ref{tm:compact-minimal-orientable} and $\tilde{\Sigma}^{(m,n)}$ its universal cover. Let $G \subseteq \mathrm{Iso}(\Sp^3_b(\kappa, \tau))$ the group of congruences of $\Sigma^{(m, n)}$ and $\tilde{G}$ the extension of G to isometries of $\tilde{\Sigma}^{(m, n)}$ by the deck
transformations of the covering. Then each element of $\tilde{G}$ extends to a congruence of the sister immersion.
\end{proposition}

\begin{proof}
Let $\phi: \tilde{\Sigma}^{(m,n)} \rightarrow E$, where $E$ is some homogeneous manifold with isometry group of dimension $4$, the sister immersion of the minimal immersion of $\tilde{\Sigma}^{(m,n)}$ in $\Sp^3_b(\kappa, \tau)$. Take a fundamental piece $\Sigma_\Gamma$. Suppose that $\tilde{\Sigma}_\Gamma$ and $\tilde{\Sigma}_\Gamma^*$ are two domain in the universal cover bijectively mapped into $\Sigma$ and $g(\Sigma^{(m,n)}_\Gamma)$ respectively, with $g \in G$. Let $\tilde{g}$ the lift of $g$ to the universal cover.
\[
\xymatrix{
\Sigma_\Gamma	\ar[d]_g & \ar[l]_{\pi} \tilde{\Sigma}_\Gamma \ar[r]^\phi \ar[d]_{\tilde{g}} \ar[dr]^\psi	& E \ar[d]^F\\
g(\Sigma_\Gamma)	& \tilde{\Sigma}_\Gamma^* \ar[l]^\pi \ar[r]_\phi & E
}
\]
Then we get two different sister immersions of $\tilde{\Sigma}_\Gamma$, the initial $\phi$ and $\psi = \phi \circ \tilde{g}$, of the same surface $\tilde{\Sigma}_\Gamma$. By the uniqueness, there exist an isometry $F: E \rightarrow E$ such that $\psi = F \circ \phi$. We claim that $F$ is a symmetry of the whole immersion $\phi: \Sigma^{(m,n)} \rightarrow E$. In fact $\phi(\Sigma^{(m,n)})$ and $(F \circ \phi)(\Sigma^{(m,n)})$ are two constant mean curvature surfaces in $E$ such that they share a open piece, $\phi(\tilde{\Sigma}_\Gamma)$, so they coincide, that is $F$ is a congruence of the immersion $\phi$.
\end{proof}

In the product space case, as we are going to see, it is possible to precise which kind of symmetry it is produced in the sister immersion. From now on we will denote $\Sp^2\times \R$ and $\h^2 \times \R$ as $M^2(\epsilon) \times \R$ if $\epsilon = 1$ and $\epsilon = -1$ respectively. Via the Daniel correspondence the shape operator $S^*$ of the constant mean curvature immersion in $M^2(\epsilon) \times \R$ is related with the shape operator $S$ of the minimal immersion in the Berger sphere as follows
\begin{equation}\label{eq:shape-operator-sister}
S^* = -JS + \tau \mathrm{Id}, \qquad \text{where } 
J = \begin{pmatrix}
0	&	-1	\\
1	&	0
\end{pmatrix}
\end{equation}
Besides, if we write the vertical Killing field as $\xi = T + \nu N$, where $N$ is a unit normal to the surface $\Sigma$ in the Berger sphere, then $(0, 1) \in T (M^2(\epsilon)\times \R) \equiv TM^2(\epsilon) \times \R$ is given by
\begin{equation}\label{eq:killing-vector-sister}
(0, 1) = -JT + \nu N^*
\end{equation}
where $N^*$ is unit normal to $\Sigma$ in $M^2(\epsilon) \times \R$.

\begin{lemma}\label{lm:simetrias-productos}
Let $\Sigma$ be a simply connected surface. Suppose that there exists an isometric minimal immersion of $\Sigma$ in $\Sp^3_b(\kappa, \tau)$. Then there is a correspondence in the sister immersions between the following symmetries of the surface $\Sigma$:
\begin{enumerate}[(a)]
	\item a reflection around a horizontal geodesic of $\Sp^3_b(\kappa, \tau)$ and
 	\item a reflection with respect a vertical plane of $M^2(\epsilon) \times \R$.
\end{enumerate}
\end{lemma}

\begin{remark}
\begin{enumerate}
	\item The relation~\eqref{eq:shape-operator-sister} remember the analogous one between constant mean curvature one surfaces in $\R^3$, considered as conjugate cousins of minimal surfaces in $\Sp^3$. In fact, this relation is the key to prove the lemma.
	
	\item Unfortunately, the relation between the shape operators in the Daniel correspondence between minimal immersions in the Berger spheres and constant mean curvature immersions in the Heisenberg group is not as good as in the correspondence to the product spaces. This behavior does not allow us to describe which kind of symmetries has the sister immersion in the Heisenberg group.
\end{enumerate}
\end{remark}

\begin{proof} 
We are going to see that (a)$\Rightarrow$(b) and the same argument goes backward. We follow the above notation and we consider $\gamma = (\beta, h) \subseteq \Sigma \subseteq M^2(\epsilon) \times \R$ parametrized by arc length and $M^2(\epsilon) \times \R \subseteq \R^3 \times \R$ (or $\R^3_1 \times \R$ if $\epsilon = -1$) isometrically immersed with unit normal along $\gamma$ given by $(\beta, 0)$. Then, we can view $J\gamma'$, where $J$ is the rotation of $\pi/2$ in the tangent plane, as a curve in $\R^4$ (or $\R^4_1$). We claim that $J\gamma'$ is constant. Firstly we get that
\[
\begin{split}
\prodesc{(J\gamma')'}{\gamma'} = -\prodesc{J\gamma'}{\nabla^\Sigma_{\gamma'}\gamma'} &= 0, \qquad \text{as $\gamma$ is a geodesic of $\Sigma$} \\
\prodesc{(J\gamma')'}{J\gamma'} &= 0, \qquad \text{as $J\gamma'$ has length $1$} \\
\prodesc{(J\gamma')'}{N^*} = -\prodesc{J\gamma'}{\df N^*(\gamma')} &= \prodesc{J\gamma'}{S^*\gamma'} =   \prodesc{J\gamma'}{-S\gamma' + \tau \gamma'} = \\
-\prodesc{\gamma'}{S\gamma'} &= 0, \qquad \text{as $\gamma$ is an asymptote line of $\Sigma$} \\
\end{split}
\]
where we take into account the relation~\eqref{eq:shape-operator-sister}. So $(J\gamma')'$ is proportional to $(\beta, 0)$ as it does not have components in the reference $\{\gamma', J\gamma', N^*\}$ of $T_{\gamma}M^2(\epsilon) \times \R$. But, because $\gamma$ is horizontal as a curve of the Berger sphere we know that $\prodesc{\gamma'}{\xi} = \prodesc{\gamma'}{T} = 0$ an so it must be $\gamma' = JT/\sqrt{1-\nu^2}$. Therefore, using~\eqref{eq:killing-vector-sister},
\[
\prodesc{\gamma'}{(0, 1)} = \frac{1}{\sqrt{1-\nu^2}}\prodesc{JT}{-JT + \nu N^*} = -\sqrt{1-\nu^2},
\]
and so $(0, 1) = -\sqrt{1 - \nu^2}\gamma' + \nu N^*$. The last equation implies that $\prodesc{J\gamma'}{(0, 1)} = 0$. Hence
\[
\prodesc{(J\gamma')'}{(\beta, 0)} = -\prodesc{J\gamma'}{(\beta', 0)} = h'\prodesc{J\gamma'}{(0, 1)} = 0
\]
where we have used that $0 = \prodesc{J\gamma'}{\gamma'} = \prodesc{J\gamma'}{(\beta', 0)} + \prodesc{J\gamma'}{h'(0,1)}$, and the claim is proved.

Finally $J\gamma' = (v, 0) \in \R^3 \times \R$. On one hand
\[
\prodesc{\gamma}{(v,0)}' = \prodesc{\gamma'}{(v,0)} = \prodesc{\gamma'}{J\gamma'} = 0
\]
so $\prodesc{\gamma}{v}$ is constant. On the other hand 
\[
\prodesc{\gamma}{(v,0)} = \prodesc{\beta}{v} = 0
\]
as $\beta$ is normal to $M^2(\epsilon)$ and $v$ is tangent. All this information says that $\gamma \subseteq P$, where $P$ is the vertical plane defined by $P = \{(p, t) \in M^2(\epsilon) \times \R:\, \prodesc{p}{v} = 0\}$. Moreover, because the tangent plane of $\Sigma$ along $\gamma$ is spanned by $\{\gamma', J\gamma' = (v,0)\}$ the surface $\Sigma$ is orthogonal to the vertical plane $P$.

Finally, as the reflexion with respect to $P$ is an isometry of $M^2(\epsilon) \times \R$ and the constant mean curvature surface $\Sigma$ is orthogonal to $P$ we can extend it by this isometry.
\end{proof}

\end{document}